# Exploring spatial nonlinearity using additive approximation

ZUDI LU[1,2,6], ARVID LUNDERVOLD[3], DAG TJøSTHEIM[4] and QIWEI YAO[1,5]

[1] *Department of Statistics, London School of Economics, London WC2A 2AE, UK*
[2] *Institute of Systems Science, Academy of Mathematics and Systems Sciences Chinese Academy of Sciences, Beijing 100080, China*
[3] *Department of Physiology, University of Bergen, 5007 Bergen, Norway*
[4] *Department of Mathematics, University of Bergen, 5007 Bergen, Norway. E-mail: dagt@mi.uib.no*
[5] *Guanghua School of Management, Peking University, Beijing 100871, China*
[6] *Department of Mathematics and Statistics, Curtin University of Technology, Perth, Australia*

We propose to approximate the conditional expectation of a spatial random variable given its nearest-neighbour observations by an additive function. The setting is meaningful in practice and requires no unilateral ordering. It is capable of catching nonlinear features in spatial data and exploring local dependence structures. Our approach is different from both Markov field methods and disjunctive kriging. The asymptotic properties of the additive estimators have been established for $\alpha$-mixing spatial processes by extending the theory of the backfitting procedure to the spatial case. This facilitates the confidence intervals for the component functions, although the asymptotic biases have to be estimated via (wild) bootstrap. Simulation results are reported. Applications to real data illustrate that the improvement in describing the data over the auto-normal scheme is significant when nonlinearity or non-Gaussianity is pronounced.

## 1. Introduction

Markov random fields and kriging are two important tools for investigating continuous spatial data. The former, including the auto-normal scheme of Besag [1] and the framework for exponential distribution families of Cressie ([5], Chapters 6 and 7), is for data on a lattice (or a graph). The latter is for irregularly positioned data; see Rivoirard [23] and Chiles and Delfiner ([4], Chapter 6). They both rely on parametric assumptions on the underlying processes. In contrast, nonparametric techniques have only found limited use in spatial modelling. This is largely due to the difficulties associated with the 'curse of dimensionality'. For example, a purely nonparametric estimation of the conditional mean at one location, given its (regularly spaced) four nearest-neighbour observations, involves four-dimensional smoothing. Although semiparametric and nonparametric autoregressive models with additive noise have proved to be successful in modelling time series, such a structure has not been available for spatial processes simply because there exists no *natural* unilateral order over a plane. On the other hand, the Markov assumption is more restrictive for spatial processes; for instance, a conditional Gaussian Markov model essentially implies linearity (Gao *et al.* [10]).

In this paper, we propose to *approximate* the conditional expectation of $Y(\mathbf{s})$, the value of a spatial process at location $\mathbf{s}$, given its $d$ nearest-neighbour observations, by an additive function, and we estimate this additive approximation by adapting the backfitting algorithm of Mammen *et al.* [15] which involves up to two-dimensional smoothing only, regardless of the value of $d$. Our approach is linked to a lattice setting. Note that data on a regular grid and measured on a continuous scale are becoming more and more common with the increasing use of computer technology. We refer to Section 5 for some tentative ideas on extending the approach to handle irregularly spaced data.

Our additive approximation may be viewed as a projection of the conditional expectation into the Hilbert space spanned by additive functions. In fact the projection principle itself does not require a lattice framework. In the context of spatial modelling, it has been used in the form of disjunctive kriging (Matheron [18, 19]; Rivoirard [23]; Chiles and Delfiner [4], Chapter 6). Disjunctive kriging projects $Y(\mathbf{s})$ into an additive space spanned by $Y(\mathbf{s}_i)$ for all $\mathbf{s}_i \neq \mathbf{s}$. Very often what is of interest is a functional of $Y(\mathbf{s})$ rather than $Y(\mathbf{s})$ itself. Nevertheless the projection principle still applies with $Y(\mathbf{s})$ replaced by $f(Y(\mathbf{s}))$ for some function $f$. For non-regularly spaced sites it is difficult to use nonparametric estimation because of the lack of repeatability of the spatial pattern of neighbours as one moves from one site to another. Instead, disjunctive kriging introduces parametric assumptions on the bivariate distributions for all pairs $(Y(\mathbf{s}_i), Y(\mathbf{s}_j))$, which then, building on appropriate *isofactorial* models, implies a parametric form for the projection of interest; see Chiles and Delfiner ([4], Chapter 6) and Rivoirard [23].

Our approach is nonparametric and pragmatic; we do not impose any explicit form on the underlying process. Instead we seek the best additive approximation to the unknown conditional expectation, which itself may not be additive. This enables us to describe local spatial dependence structure with a potential application to texture analysis. For example, the nonlinear structure demonstrated in the additive approximation for the straw data in Section 4.2 is beyond the capacity of an auto-normal fitting and would be difficult to describe using disjunctive kriging. Our approach also provides a vehicle for testing isotropy and/or linearity; see the bootstrap test in Example 2 in Section 4.1. This may serve as a guide for choosing a parametric model. Of course, those advantages come at some cost. For example, abandoning the Markov framework implies that Markov chain Monte Carlo and other important analytical tools are not at our disposal. This may be a severe obstacle when dealing with non-stationary spatial processes.

Another way of circumventing the curse of dimensionality is to use semiparametric (partially linear) additive approximation if some components are found to be linear; this is explored by Gao *et al.* [10] with the marginal integration technique (Linton and Nielsen [14]; Newey [20]; Tjøstheim and Auestad [24]). The marginal integration method is less efficient in practice than back-fitting when $d$ is large, in spite of its good asymptotic properties (Fan *et al.* [9]). Both methods require density estimation. It should be noted that nonparametric density estimation for spatial processes can be traced back at least to Diggle [7] and Diggle and Marron [8]. For more recent developments, see Carbon *et al.* [3], Hallin *et al.* [11, 12, 13] and Yao [26].

The rest of the paper is organized as follows. The methodology is laid out in Section 2. Asymptotic properties are stated in Section 3. The asymptotic distributions of the estimators are used to construct pointwise confidence intervals for component functions in the additive approximation, although the asymptotic biases are estimated via wild bootstrap. Numerical illustrations with both simulated and real data sets are reported in Section 4. A brief discussion on possible extension to handling irregularly spaced data is presented in Section 5. All technical proofs are relegated to the Appendix.

# 2. Methodology

## 2.1. Least-squares additive approximation

Suppose $\{Y(\mathbf{s})\}$ is a strictly stationary process defined on a two-dimensional lattice, that is, $\mathbf{s} \equiv (u,v) \in \mathbb{Z}^2$, where $\mathbb{Z}$ denotes the set consisting of all integers. Let $\mathbf{i}_1, \ldots, \mathbf{i}_d$ be $d$ fixed neighbour points of $(0,0)$ in $\mathbb{Z}^2$, and $\mathbf{x} = (x_1, \ldots, x_d)^\top \in \mathbb{R}^d$. It is of interest to approximate the conditional expectation

$$m(\mathbf{x}) \equiv \mathrm{E}\{Y(\mathbf{s}) | Y(\mathbf{s} - \mathbf{i}_\ell) = x_\ell, \ell = 1, \ldots, d\} \tag{2.1}$$

by an additive form

$$m_0 + m_1(x_1) + \cdots + m_d(x_d). \tag{2.2}$$

We seek the optimal approximation in a least-squares sense; see (2.4) below. To make the terms in (2.2) identifiable, we require $\int m_j(y) f_0(y)\,\mathrm{d}y = 0$, $j = 1, \ldots, d$, where $f_0(\cdot)$ denotes the marginal density function of $Y(\mathbf{s})$. If $m(\cdot)$ itself is of the form (2.2), it is easy to see that

$$m_j(x_j) = \mathrm{E}\{Y(\mathbf{s}) | Y(\mathbf{s} - \mathbf{i}_j) = x_j\} - m_0 - \sum_{\substack{1 \le \ell \le d \\ \ell \ne j}} \mathrm{E}[m_\ell\{Y(\mathbf{s} - \mathbf{i}_\ell)\} | Y(\mathbf{s} - \mathbf{i}_j) = x_j]. \tag{2.3}$$

In general, we obtain an optimum additive approximation by minimizing

$$E\left[Y(\mathbf{s}) - m_0 - \sum_{\ell=1}^{d} m_\ell\{Y(\mathbf{s} - \mathbf{i}_\ell)\}\right]^2, \tag{2.4}$$

or equivalently,

$$\mathrm{E}\left[m\{\mathbf{X}(\mathbf{s})\} - m_0 - \sum_{\ell=1}^{d} m_\ell\{Y(\mathbf{s} - \mathbf{i}_\ell)\}\right]^2, \tag{2.5}$$

over $m_0 + \sum_{\ell=1}^{d} m_\ell(\cdot) \in \mathcal{F}_{\text{add}}$, where $\mathbf{X}(\mathbf{s}) = \{Y(\mathbf{s}-\mathbf{i}_1), \ldots, Y(\mathbf{s}-\mathbf{i}_d)\}^\top$ and

$$\mathcal{F}_{\text{add}} = \left\{ m_0 + \sum_{\ell=1}^{d} m_\ell(x_\ell) \middle| m_0 \in \mathbb{R}, \int m_\ell(y) f_0(y)\, \mathrm{d}y = 0 \text{ for } 1 \le \ell \le d \right\}. \tag{2.6}$$

## 2.2. Estimators

We now spell out how to estimate the best additive approximation for the conditional expectation (2.1). To simplify notation, we assume that observations $\{(Y(\mathbf{s}_\ell), \mathbf{X}(\mathbf{s}_\ell)), 1 \le \ell \le N\}$ are available. Furthermore, we assume that those data are taken from a rectangle in $\mathbb{Z}^2$, for example,

$$\{\mathbf{s}_1, \ldots, \mathbf{s}_N\} = \{(u, v) : u = 1, \ldots, N_1, v = 1, \ldots, N_2\}, \tag{2.7}$$

where $N_1 N_2 = N$. Other sampling schemes are possible; see the remark at the end of Section 3.1 below.

In practice we replace (2.5) by

$$\int \left\{ \widehat{m}(\mathbf{x}) - m_0 - \sum_{\ell=1}^{d} m_\ell(x_\ell) \right\}^2 \widehat{f}(\mathbf{x})\, \mathrm{d}\mathbf{x}, \tag{2.8}$$

where $\mathbf{x} = (x_1, \ldots, x_d)^\top$, and

$$\widehat{f}(\mathbf{x}) = \frac{1}{N} \sum_{\ell=1}^{N} K_h\{\mathbf{x} - \mathbf{X}(\mathbf{s}_\ell)\}, \qquad K_h(\mathbf{x}) = \prod_{j=1}^{d} K_h(x_j),$$

$$\widehat{m}(\mathbf{x}) = \frac{\widehat{r}(\mathbf{x})}{\widehat{f}(\mathbf{x})}, \qquad \widehat{r}(\mathbf{x}) = \frac{1}{N} \sum_{\ell=1}^{N} Y(\mathbf{s}_\ell) K_h\{\mathbf{x} - \mathbf{X}(\mathbf{s}_\ell)\}. \tag{2.9}$$

In the above expression, $K_h(x) = h^{-1}K(x/h)$, $K(\cdot)$ is a density function on $\mathbb{R}$, and $h > 0$ is the bandwidth. We also define

$$\widehat{f}_j(x_j) = \frac{1}{N} \sum_{\ell=1}^{N} K_h\{x_j - Y(\mathbf{s}_\ell - \mathbf{i}_j)\}, \qquad \widehat{m}_j(x_j) = \frac{1}{N \widehat{f}_j(x_j)} \sum_{\ell=1}^{N} Y(\mathbf{s}_\ell) K_h\{x_j - Y(\mathbf{s}_\ell - \mathbf{i}_j)\},$$

$$\widehat{f}_{jk}(x_j, x_k) = \frac{1}{N} \sum_{\ell=1}^{N} K_h\{x_j - Y(\mathbf{s}_\ell - \mathbf{i}_j)\} K_h\{x_k - Y(\mathbf{s}_\ell - \mathbf{i}_k)\}.$$

Note that $\widehat{f}_j$ and $\widehat{f}_{jk}$ are the marginal density functions from the joint density $\widehat{f}$.

Let $\{\widetilde{m}_l\} \in \widehat{\mathcal{F}}_{\text{add}}$ be a minimizer of (2.8), where

$$\widehat{\mathcal{F}}_{\text{add}} = \left\{ m_0 + \sum_{\ell=1}^{d} m_\ell(x_\ell) \middle| m_0 \in \mathbb{R}, \int m_\ell(y) \widehat{f}_\ell(y)\, \mathrm{d}y = 0 \text{ for } 1 \le \ell \le d \right\}.$$

Then the least-squares property implies that

$$\int \left\{ \widehat{m}(\mathbf{x}) - \widetilde{m}_0 - \sum_{\ell=1}^{d} \widetilde{m}_\ell(x_\ell) \right\} \{\widetilde{m}_j(x_j) - m_j(x_j)\} \widehat{f}(\mathbf{x})\, \mathrm{d}\mathbf{x} = 0$$

for any $m_j(\cdot)$, $j = 0, 1, \ldots, d$. (We write $m_0(\cdot) \equiv m_0$.)

This leads to

$$\widetilde{m}_0 = \int \left\{ \widehat{m}(\mathbf{x}) - \sum_{\ell=1}^{d} \widetilde{m}_\ell(x_\ell) \right\} \widehat{f}(\mathbf{x})\, \mathrm{d}\mathbf{x} = \int \widehat{m}(\mathbf{x}) \widehat{f}(\mathbf{x})\, \mathrm{d}\mathbf{x} = \frac{1}{N} \sum_{\ell=1}^{N} Y(\mathbf{s}_\ell) \equiv \bar{Y}, \tag{2.10}$$

and for $j = 1, \ldots, d$,

$$\widetilde{m}_j(x_j) = \widehat{m}_j(x_j) - \widetilde{m}_0 - \sum_{\ell \neq j} \int \widetilde{m}_\ell(x_\ell) \frac{\widehat{f}_{j\ell}(x_j, x_\ell)}{\widehat{f}_j(x_j)}\, \mathrm{d}x_\ell. \tag{2.11}$$

We always follow the convention that $x/y$ equals 0 if $y = 0$. It is easy to see that (2.11) intimately resembles (2.3). It also naturally leads to the following backfitting algorithm: in the $j$th step of the $r$th iteration cycle we define

$$\begin{aligned} \widetilde{m}_j^{(r)}(x_j) = \widehat{m}_j(x_j) - \bar{Y} &- \sum_{\ell < j} \int \widetilde{m}_\ell^{(r)}(x_\ell) \frac{\widehat{f}_{j\ell}(x_j, x_\ell)}{\widehat{f}_j(x_j)}\, \mathrm{d}x_\ell \\ &- \sum_{\ell > j} \int \widetilde{m}_\ell^{(r-1)}(x_\ell) \frac{\widehat{f}_{j\ell}(x_j, x_\ell)}{\widehat{f}_j(x_j)}\, \mathrm{d}x_\ell. \end{aligned} \tag{2.12}$$

We choose Nadaraya–Watson (i.e. local constant) estimation to keep our exposition as simple as possible. For general discussion of smoothing backfitting algorithms, including the one based on more efficient local linear estimation, we refer to Mammen *et al.* [15] and Nielsen and Sperlich [21]. More recently, Mammen and Park [17] showed that if in (2.12) $\widehat{m}_j$ is replaced by a marginal local linear estimator, and $\widehat{f}_{j\ell}/\widehat{f}_j$ is replaced by a more sophisticated functional constructed using a convolution kernel, the resulting backfitting method is asymptotically as efficient as the one based on local linear estimation.

Finally, we remark that the minimizer of (2.8) is the same as the minimizer of

$$\frac{1}{N} \int \sum_{j=1}^{N} \left\{ Y(\mathbf{s}_j) - m_0 - \sum_{\ell=1}^{d} m_\ell(x_\ell) \right\}^2 K_h\{\mathbf{x} - \mathbf{X}(\mathbf{s}_j)\}\, \mathrm{d}\mathbf{x}.$$

This lends support to the use of a simple leave-one-out cross-validation bandwidth estimator:

$$\widehat{h} = \arg\min_h \sum_{j=1}^{N} \left[ Y(\mathbf{s}_j) - \widetilde{m}_{0,-j} - \sum_{\ell=1}^{d} \widetilde{m}_{\ell,-j}\{Y(\mathbf{s}_j - \mathbf{i}_\ell)\} \right]^2, \tag{2.13}$$

where $\widetilde{m}_{\ell,-j}$ is the backfitting estimator of $m_\ell$ without the $j$th observation $(Y(\mathbf{s}_j), \mathbf{X}(\mathbf{s}_j))$. Nielsen and Sperlich [21] proposed some modifications to make this bandwidth selector computationally more efficient. Three other data-driven bandwidth selectors for additive modelling based on backfitting were proposed in Mammen and Park [16].

# 3. Asymptotic properties

## 3.1. Regularity conditions

In order to present asymptotic results, we define the $\alpha$-mixing coefficients for spatial processes first. For any $A \subset \mathbb{Z}^2$, let $\mathcal{F}(A)$ denote the $\sigma$-algebra generated by $\{(\mathbf{X}(\mathbf{s}), Y(\mathbf{s})), \mathbf{s} \in A\}$. We write $|A|$ for the number of elements in $A$. For any $A, B \subset \mathbb{Z}^2$, define

$$\alpha(A,B) = \sup_{U \in \mathcal{F}(A), V \in \mathcal{F}(B)} |P(UV) - P(U)P(V)|,$$

and $d(A,B) = \min\{\|\mathbf{s}_1 - \mathbf{s}_2\| \,|\, \mathbf{s}_1 \in A, \mathbf{s}_2 \in B\}$, where $\|\cdot\|$ denotes the Euclidean norm. We may define an $\alpha$-mixing coefficient for the process $\{(\mathbf{X}(\mathbf{s}), Y(\mathbf{s}))\}$ as

$$\alpha(k;i,j) = \sup_{A,B \subset \mathbb{Z}^2} \{\alpha(A,B) \,|\, |A| \le i, |B| \le j, d(A,B) \ge k\}, \tag{3.1}$$

where $i, j, k$ are positive integers and $i, j$ may be infinite. For further discussions on the mixing for spatial processes, we refer to Section 1.3.1 of Doukhan [6] and Section 2.1 of Yao [26] and references within.

Let $C$ denote some positive generic constant which may be different at different places. The following regularity conditions are imposed.

(C1) The density functions $f$ of $Y(\mathbf{s})$ and $f_{jk}$ of $\{Y(\mathbf{s} - \mathbf{i}_j), Y(\mathbf{s} - \mathbf{i}_k)\}$ have continuous second derivatives, and are bounded from above by a constant independent of $\mathbf{i}_j - \mathbf{i}_k$. The conditional expectation $m_j(\cdot)$ has continuous first derivative. The density functions of $\mathbf{X}(\mathbf{s})$ conditional on $Y(\mathbf{s})$, and $\{\mathbf{X}(\mathbf{i}), \mathbf{X}(\mathbf{s})\}$ conditional on $\{Y(\mathbf{i}), Y(\mathbf{s})\}$ are bounded from above. Furthermore, for some $\lambda > 0$ and $N^{-\lambda+3/2}h^{-1/2} \to 0$,

$$\mathrm{E}\{\exp(\lambda|Y(\mathbf{s})|)\} < \infty. \tag{3.2}$$

(C2) The kernel function $K(\cdot)$ is symmetric, compactly supported and Lipschitz continuous.

(C3) As $N = N_1 N_2 \to \infty$, $N_1/N_2$ converges to a positive and finite constant, the bandwidth $h \to 0$ and

$$N^{\beta-5}h^{\beta+5}(\log N)^{-(3\beta+7)} \to \infty, \tag{3.3}$$

where $\beta > 5$ is a constant.

(C4) $\alpha(k;k',\infty) \le Ck^{-\beta}$ for any $k$ and $k' = O(k^2)$. Furthermore, $\sum_{k=1}^{\infty} k^{d-1}\alpha(k;j,\ell) < \infty$ for some $j+\ell \le 4$, $\alpha(k;1,\infty) = o(k^{-d})$ and $\sum_{k=1}^{\infty} k^{d-1}\alpha(k;1,1)^{(\delta-2)/\delta} < \infty$ for some $\delta > 2$.

Conditions (C1)–(C2) are standard in kernel estimation. Both the assumption of the existence of the moment generating function of $|Y(\mathbf{s})|$ and (3.3), which imply the optimum uniform convergence rates (A.1) and (A.2), can be relaxed at the cost of lengthy arguments. On the other hand, for causal and invertible (under the half-plane order) spatial ARMA processes satisfying some mild conditions, $\alpha(k;k',\infty)$ decays at an exponential rate as $k \to \infty$ (Remark 2.1 of Yao [26]). Therefore, condition (C4) is fulfilled. For optimum bandwidth $h = O(N^{-1/5})$, (3.3) is fulfilled for $\beta > 7.5$. Condition (C3) requires two sides of the sampling rectangle to increase to infinity. In fact this assumption can be relaxed. For example, our theoretical results will still hold if the observations were taken over a connected region in $\mathbb{Z}^2$, and the ratio of the minimal side length of the squares containing the region to the maximal side length of the squares contained in the region converges to a constant in the interval $(0,\infty)$. For a general discussion on the condition of sampling sets, we refer to Perera [22].

## 3.2. Convergence of backfitting

Backfitting techniques have proved effective in handling complex model fitting. However, its convergence is typically difficult to handle. We apply Theorem 1 of Mammen *et al.* [15] to show that a modified version of backfitting (2.12) converges. The modification is in line with Section 5 of Mammen *et al.* [15] in order to fulfil certain regularity conditions which simplify technical arguments substantially.

Let $A \subset \mathbb{R}^d$ be a compact set on which the density function $f_{1,\dots,d}(\cdot)$ of $\mathbf{X}(\mathbf{s})$ is positive. Define

$$p(\mathbf{x}) \equiv p_{1,\dots,d}(\mathbf{x}) = \frac{f_{1,\dots,d}(\mathbf{x}) I(\mathbf{x} \in A)}{P\{\mathbf{X}(\mathbf{s}) \in A\}}.$$

Then $p(\cdot)$ is a density function on $\mathbb{R}^d$. As an illustration, Mammen *et al.* [15] chose $A = [0,1]^d$. Since the components of $\mathbf{X}(\mathbf{s})$ are dependent, sets of cylinder type are not always relevant. (For example, the support of $(X_t, X_{t-1})$ for linear AR(1) time series would be around a line segment with non-zero slope.) Denote by $p_j(x_j)$ and $p_{jk}(x_j,x_k)$ respectively the $j$th univariate and the $(j,k)$th bivariate marginal density functions of $p(\mathbf{x})$. We require the following consistency condition on the set $A$.

(C5) There exist compact sets $A_j \subset \{f(x_j) > 0\}$ and $A_{jk} \subset \{f_{jk}(x_j,x_k) > 0\}$ such that for $1 \le j,k \le d$ and $j \ne k$,

$$p_j(x_j) = \frac{f(x_j) I(x_j \in A_j)}{P\{Y(\mathbf{s}) \in A_j\}}, \qquad p_{jk}(x_j,x_k) = \frac{f_{jk}(x_j,x_k) I\{(x_j,x_k) \in A_{jk}\}}{P[\{Y(\mathbf{s}-\mathbf{i}_j), Y(\mathbf{s}-\mathbf{i}_k)\} \in A_{jk}]}.$$

Due to stationarity, a relevant set $A$ often exhibits certain symmetries. For example, we may observe $A_i = A_j$ and $p_i(\cdot) = p_j(\cdot)$ for all $i$ and $j$.

Differently from Section 5 of Mammen *et al.* [15], we define estimators for $p_j$ and $p_{jk}$ as follows:

$$\widehat{p}_j(x_j) = I(x_j \in A_j)\frac{\sum_{\ell=1}^N K_h\{x_j - Y(\mathbf{s}_\ell - \mathbf{i}_j)\}}{\sum_{\ell=1}^N I\{Y(\mathbf{s}_\ell - \mathbf{i}_j) \in A_j\}}, \tag{3.4}$$

$$\widehat{p}_{jk}(x_j, x_k) = I\{(x_j, x_k) \in A_{jk}\}\frac{\sum_{\ell=1}^N K_h\{x_j - Y(\mathbf{s}_\ell - \mathbf{i}_j)\}K_h\{x_k - Y(\mathbf{s}_\ell - \mathbf{i}_k)\}}{\sum_{\ell=1}^N I[\{Y(\mathbf{s}_\ell - \mathbf{i}_j), Y(\mathbf{s}_\ell - \mathbf{i}_k)\} \in A_{jk}]}. \tag{3.5}$$

Obviously $\widehat{p}_j$ and $\widehat{p}_{jk}$ are consistent estimators for $p_j$ and $p_{jk}$, respectively. Note that $K(\cdot)$ is compactly supported. For $x_j \in A_j$, $K_h\{x_j - Y(\mathbf{s}_\ell - \mathbf{i}_j)\}$ may be non-zero for sufficiently large $N$ only if $Y(\mathbf{s}_\ell - \mathbf{i}_j) \in A_j'$, where $A_j'$ is a compact set sandwiched between $A_j$ and $\{f(x_j) > 0\}$. Therefore, similarly to Mammen *et al.* [15], we effectively only use the observations in a compact set when estimating $\widehat{p}_j$. It is possible now that $\int \widehat{p}_{jk}(x_j, x_k)\,\mathrm{d}x_k \neq \widehat{p}_j(x_j)$. Similarly to Mammen *et al.* [15], we modify the backfitting procedure (2.11) and (2.12) accordingly:

$$\widetilde{m}_j(x_j) = \widehat{m}_j(x_j) - \widetilde{m}_{0,j} - \sum_{\ell \neq j}\int \widetilde{m}_\ell(x_\ell)\left\{\frac{\widehat{p}_{j\ell}(x_j, x_\ell)}{\widehat{p}_j(x_j)} - \frac{\int \widehat{p}_{j\ell}(u, x_\ell)\,\mathrm{d}u}{\int \widehat{p}_j(u)\,\mathrm{d}u}\right\}\mathrm{d}x_\ell, \tag{3.6}$$

$$\begin{aligned}\widetilde{m}_j^{(r)}(x_j) = \widehat{m}_j(x_j) - \widetilde{m}_{0,j} &- \sum_{\ell < j}\int \widetilde{m}_\ell^{(r)}(x_\ell)\left\{\frac{\widehat{p}_{j\ell}(x_j, x_\ell)}{\widehat{p}_j(x_j)} - \frac{\int \widehat{p}_{j\ell}(u, x_\ell)\,\mathrm{d}u}{\int \widehat{p}_j(u)\,\mathrm{d}u}\right\}\mathrm{d}x_\ell \\ &+ \sum_{\ell > j}\int \widetilde{m}_\ell^{(r-1)}(x_\ell)\left\{\frac{\widehat{p}_{j\ell}(x_j, x_\ell)}{\widehat{p}_j(x_j)} - \frac{\int \widehat{p}_{j\ell}(u, x_\ell)\,\mathrm{d}u}{\int \widehat{p}_j(u)\,\mathrm{d}u}\right\}\mathrm{d}x_\ell,\end{aligned} \tag{3.7}$$

where $\widetilde{m}_{0,j} = \int \widehat{m}_j(x)\widehat{p}_j(x)\,\mathrm{d}x / \int \widehat{p}_j(y)\,\mathrm{d}y$. Note that, for $x_j \in A_j$, $\widehat{m}_j(x_j)$ defined in (2.9) may be written as

$$\widehat{m}_j(x_j) = I(x_j \in A_j)\frac{\sum_{\ell=1}^N Y(\mathbf{s}_\ell)K_h\{x_j - Y(\mathbf{s}_\ell - \mathbf{i}_j)\}}{\widehat{p}_j(x_j)\sum_{\ell=1}^N I\{Y(\mathbf{s}_\ell - \mathbf{i}_j) \in A_j\}}.$$

The following theorem indicates that this backfitting procedure converges exponentially fast.

**Theorem 1.** *Under conditions* (C1)–(C5)*, with probability tending to 1, there exists a unique solution* $\{\widetilde{m}_j\}$ *of* (3.6)*, and further, for* $r \geq 1$ *and* $\mathbf{x} = (x_1, \ldots, x_d)^\top$ *being an inner point of* $A$*,*

$$\sum_{j=1}^d \int \{\widetilde{m}_j^{(r)}(x_j) - \widetilde{m}_j(x_j)\}^2 p_j(x_j)\,\mathrm{d}x_j \leq C\rho^{2r}\left(1 + \sum_{j=1}^d \int \{\widetilde{m}_j^{(0)}(x_j)\}^2 p_j(x_j)\,\mathrm{d}x_j\right),$$

*where* $\rho \in (0,1), C > 0$ *are constants,* $\widetilde{m}_j^{(r)}(x_j)$ *is defined by* (3.7)*, and* $\{\widetilde{m}_j^{(0)}(x_j)\}$ *are the initial values of the backfitting algorithm.*

### 3.3. Asymptotic normality

Henceforth, we assume that $\mathbf{x} = (x_1, \dots, x_d)^\top$ is an inner point of $A$. Let $\varepsilon(\mathbf{s}) = Y(\mathbf{s}) - m\{\mathbf{X}(\mathbf{s})\}$, and $m^o(\mathbf{x}) = m_0^o + \sum_{j=1}^d m_j^o(x_j)$ be the minimizer of (2.5) over

$$\mathcal{F}'_{\text{add}} = \left\{ m(\mathbf{x}) = m_0 + \sum_{j=1}^d m_j(x_j) \middle| m_0 \in \mathbb{R}, \int m_j(y) p_j(y)\, \mathrm{d}y = 0 \text{ for } 1 \le j \le d \right\}.$$

Then $\{m_0^o, m_1^o(\cdot), \dots, m_d^o(\cdot)\}$ is uniquely determined by the least-squares property. Define

$$\beta(\mathbf{x}) = \sum_{j=1}^d \left\{ \dot{m}_j^o(x_j) \frac{\partial}{\partial x_j} \log p(\mathbf{x}) + \frac{1}{2} \ddot{m}_j^o(x_j) \right\} \int u^2 K(u)\, \mathrm{d}u, \tag{3.8}$$

$$\widehat{\mu}_j(x_j) = m_j^o(x_j) + \sum_{k \neq j} \int m_k^o(x_k) \frac{\widehat{p}_{jk}(x_j, x_k)}{\widehat{p}_j(x_j)} \mathrm{d}x_k + h^2 \int \beta(\mathbf{x}) \frac{p(\mathbf{x})}{p_j(x_j)} \prod_{k \neq j} \mathrm{d}x_k.$$

Let $\beta_0 + \sum_{j=1}^d \beta_j(x_j)$ be the minimizer of

$$\int \left\{ \beta(\mathbf{x}) - \beta_0 - \sum_{j=1}^d \beta_j(x_j) \right\}^2 p(\mathbf{x})\, \mathrm{d}\mathbf{x} \tag{3.9}$$

over $\mathcal{F}'_{\text{add}}$.

**Theorem 2.** *Let conditions* (C1)–(C5) *hold, and* $h = CN^{-1/5}$. *Then*

$$\sqrt{Nh} \begin{pmatrix} \widetilde{m}_1(x_1) - m_1^o(x_1) - h^2 \beta_1(x_1) \\ \vdots \\ \widetilde{m}_d(x_d) - m_d^o(x_d) - h^2 \beta_d(x_d) \end{pmatrix} \xrightarrow{D} N(0, \mathbf{\Sigma}(\mathbf{x})),$$

*where* $\mathbf{\Sigma}(\mathbf{x})$ *is a diagonal matrix with*

$$\sigma_j(x_j)^2 \equiv \operatorname{var}[Y(\mathbf{s}) - m^o\{\mathbf{X}(\mathbf{s})\} | Y(\mathbf{s} - \mathbf{i}_j) = x_j] \int K(u)^2\, \mathrm{d}u / f_j(x_j) \tag{3.10}$$

*as its* $j$*th main diagonal element.*

***Remark 1.*** (i) Although we do not assume the true conditional expectation $m(\mathbf{x})$ defined in (2.1) to be of additive form, the estimators do not have extra biases due to the discrepancy between $m(\mathbf{x})$ and its best additive approximation $m^o(\mathbf{x})$. This is due to the 'orthogonality'

$$\int \{m(\mathbf{x}) - m^o(\mathbf{x})\} p(\mathbf{x}) \prod_{k \neq j} \mathrm{d}x_k = 0, \qquad 1 \le j \le d, \tag{3.11}$$

which is guaranteed by the least-squares property that $m^o(\cdot)$ is the minimizer of (2.5) over $\mathcal{F}_{\text{add}}$. On the other hand, the variance in (3.10) is equal to

$$(\text{var}[Y(\mathbf{s}) - m\{\mathbf{X}(\mathbf{s})\}|Y(\mathbf{s}-\mathbf{i}_j) = x_j] \\ + \text{var}[m\{\mathbf{X}(\mathbf{s})\} - m^o\{\mathbf{X}(\mathbf{s})\}|Y(\mathbf{s}-\mathbf{i}_j) = x_j]) \int K(u)^2\,\mathrm{d}u \Big/ p_j(x_j).$$

The second term in the above expression disappears when $m(\mathbf{x})$ itself is an additive function.

(ii) The proof of Theorem 2 entails that

$$\widetilde{m}_j(x_j) - m_j^o(x_j) - h^2\beta_j(x_j) + o_p(h^2) \\ = \frac{1}{N\widehat{p}_j(x_j)} \sum_{\ell=1}^{N} [Y(\mathbf{s}_\ell) - m^o\{\mathbf{X}(\mathbf{s}_\ell)\}] K_h\{x_j - Y(\mathbf{s}_\ell - \mathbf{i}_j)\}. \tag{3.12}$$

By Theorem 2, an approximate 95% pointwise confidence interval for $m_j^o(x_j)$ would be of the form $\widetilde{m}_j(x_j) - h^2\beta_j(x_j) \pm 1.96\sigma_j(x_j)/\sqrt{nh}$. However, the quantities $\beta_j(x_j)$ and $\sigma_j(x_j)$ are unknown in practice. Furthermore, it is rather difficult to estimate $\beta_j(\cdot)$; see (3.9) and (3.8). We now outline a heuristic method based on wild bootstrapping to estimate the bias $\beta_j(x_j)$ and the variance $\sigma_j(x_j)^2$.

Let $\{\varepsilon(\mathbf{s})\}$ be independent and identically distributed random variables with mean 0 and variance 1. Draw a bootstrap sample $Y(\mathbf{s}_1)^*, \ldots, Y(\mathbf{s}_N)^*$ from the model

$$Y(\mathbf{s})^* = \widetilde{m}\{\mathbf{X}(\mathbf{s})\} + \varepsilon(\mathbf{s})[Y(\mathbf{s}) - \widetilde{m}\{\mathbf{X}(\mathbf{s})\}], \tag{3.13}$$

where $\widetilde{m}(\mathbf{x}) = \widetilde{m}_0 + \sum_{j=1}^{d} \widetilde{m}_j(x_j)$. Then $\mathrm{E}^*\{Y(\mathbf{s})^*|\mathbf{X}(\mathbf{s}) = \mathbf{x}\} = \widetilde{m}(\mathbf{x})$. Let $\{m_j^*\}$ be the estimators obtained in the same way as $\{\widetilde{m}_j\}$ but with sample $\{Y(\mathbf{s}_j), \mathbf{X}(\mathbf{s}_j)\}$ replaced by $\{Y(\mathbf{s}_j)^*, \mathbf{X}(\mathbf{s}_j)\}$. It may be shown that $m_j^*(x_j) - \widetilde{m}_j(x_j)$ shares the same asymptotic distribution as $\widetilde{m}_j(x_j) - m_j^o(x_j)$; see Theorem 2 above. Hence, we may use the sample mean and the sample variance of $m_j^*(x_j) - \widetilde{m}_j(x_j)$ obtained in a repeated bootstrap sampling (with a large number of replications) as the estimates for the mean and the variance of $\widetilde{m}_j(x_j) - m_j^o(x_j)$. Combining with Theorem 2, this leads to an approximate 95% pointwise confidence interval for $m_j^o(x_j)$ $(1 \le j \le d)$:

$$2\widetilde{m}_j(x_j) - \bar{m}_j^*(x_j) \pm 1.96 s_j^*(x_j), \tag{3.14}$$

where $\bar{m}_j^*(x_j)$ is the sample mean of $m_j^*(x_j)$ in the repeated bootstrap sampling, and $s_j^*(x_j)$ is the sample standard deviation of $m_j^*(x_j) - \widetilde{m}_j(x_j)$.

***Remark 2.*** Note that the conditional expectation $\mathrm{E}^*\{Y(\mathbf{s})^*|\mathbf{X}(\mathbf{s}) = \mathbf{x}\} = \widetilde{m}(\mathbf{x})$ is an additive function, while $\mathrm{E}\{Y(\mathbf{s})|\mathbf{X}(\mathbf{s}) = \mathbf{x}\} = m(\mathbf{x})$ may not be. This makes it difficult to construct confidence intervals directly based on bootstrapping. The confidence interval (3.14) is based on the asymptotic normality of the estimator $\widetilde{m}_j(x_j)$. Bootstrapping was

merely employed to estimate the unknown asymptotic bias $\beta_j(x_j)$ and variance $\sigma_j(x_j)^2$, which relied on the fact that $\beta_j(x_j)$ is determined by the best additive approximation $m^o(\mathbf{x})$ of $m(\mathbf{x})$ instead of $m(\mathbf{x})$ itself; see (3.9) and (3.8). On the other hand, it may be shown that

$$\begin{aligned}
&m_j^*(x_j) - \widetilde{m}_j(x_j) - h^2\beta_j(x_j) + o_p(h^2) \\
&= \frac{1}{N\widehat{p}_j(x_j)} \sum_{\ell=1}^{N} [Y(\mathbf{s}_\ell)^* - \widetilde{m}\{\mathbf{X}(\mathbf{s}_\ell)\}] K_h\{x_j - Y(\mathbf{s}_\ell - \mathbf{i}_j)\} \\
&= \frac{1}{N\widehat{p}_j(x_j)} \sum_{\ell=1}^{N} \varepsilon(\mathbf{s})[Y(\mathbf{s}_\ell) - \widetilde{m}\{\mathbf{X}(\mathbf{s}_\ell)\}] K_h\{x_j - Y(\mathbf{s}_\ell - \mathbf{i}_j)\} \\
&= \frac{1}{N\widehat{p}_j(x_j)} \sum_{\ell=1}^{N} \varepsilon(\mathbf{s})[Y(\mathbf{s}_\ell) - m^o\{\mathbf{X}(\mathbf{s}_\ell)\}] K_h\{x_j - Y(\mathbf{s}_\ell - \mathbf{i}_j)\}\{1 + o_p(1)\};
\end{aligned}$$

see (3.12). This would ensure that the bootstrap estimator admits the same asymptotic variance.

# 4. Numerical properties

## 4.1. Simulation

In this section, we illustrate the proposed backfitting procedure with two examples: a unilateral additive model under the half-plane order (Whittle [25]), and a (bilateral) auto-normal model (Besag [1]). The bandwidth selection procedure (2.13) is implemented in Example 1. In Example 2 we apply a parametric bootstrap test to test the null hypothesis of an auto-normal model. In the numerical examples let $K$ be a Gaussian kernel.

***Example 1*** *Unilateral additive model.* Consider the model

$$Y(u,v) = \sin\{Y(u-1,v)\} + \cos\{Y(u,v-1)\} + e(u,v), \tag{4.1}$$

where $e(u,v)$ are independent $N(0,1)$ random variables. Hence

$$\begin{aligned}
&\mathrm{E}\{Y(u,v) | Y(u-1,v), Y(u,v-1), Y(u-1,v-1)\} \\
&\quad = m_0 + m_1\{Y(u-1,v)\} + m_2\{Y(u,v-1)\} + m_3\{Y(u-1,v-1)\}
\end{aligned} \tag{4.2}$$

with $m_0 = E\{Y(u,v)\}$, $m_3(\cdot) \equiv 0$, and

$$m_1(x) = \sin(x) - \mathrm{E}[\sin\{Y(u,v)\}], \qquad m_2(x) = \cos(x) - \mathrm{E}[\cos\{Y(u,v)\}].$$

We drew 100 samples from model (4.1) on the rectangle $\{(u,v) : 1 \le u \le 24, 1 \le v \le 28\}$. For each sample we estimated the component functions $m_j(\cdot)$ for $j = 1,2,3$ with the bandwidths $h$ chosen automatically by the leave-one-out procedure (2.13). The boxplots of the estimated curves over 13 regular grid points are presented in Figure 1. While the estimation is accurate overall, the variation of the estimation is larger at the both ends due to boundary effects. The mean and variance of the selected bandwidths over 100 replications are 0.416 and 0.064, respectively.

***Example 2*** *Besag's first-order auto-normal scheme.* Let the conditional distribution of $Y(u,v)$ given $\{Y(i,j),(i,j)\neq(u,v)\}$ be normal with mean

$$\begin{aligned} &\mathrm{E}\{Y(u,v)|Y(i,j),(i,j)\neq(u,v)\} \\ &\quad=\theta_1\{Y(u-1,v)+Y(u+1,v)\}+\theta_2\{Y(u,v-1)+Y(u,v+1)\} \end{aligned} \tag{4.3}$$

and variance $\mathrm{var}\{Y(u,v)|Y(i,j),(i,j)\neq(u,v)\}=1$, where $\theta_1=0.2$ and $\theta_2=0.25$. Now the $\{Y(u,v)\}$ are jointly normal with common mean 0, and

$$\begin{aligned} &\mathrm{E}\{Y(u,v)|Y(u-1,v)=x_1,Y(u,v-1)=x_2,Y(u+1,v)=x_3,Y(u,v+1)=x_4\} \\ &\quad=m_1(x_1)+m_2(x_2)+m_3(x_3)+m_4(x_4) \end{aligned}$$

with $m_1(x)=m_3(x)=\theta_1 x$, $m_2(x)=m_4(x)=\theta_2 x$; see Besag [1].

We conducted a simulation with 500 replications. For each sample taken on the rectangle $\{(u,v):1\leq u,v\leq 20\}$, we applied the backfitting algorithm to estimate $m_j(\cdot)$. The boxplots of the estimated curves over 11 grid points are presented in Figure 2. To speed up the computation, we used a fixed bandwidth $h=0.4$. The linearity of $m_j(\cdot)$ is evident in Figure 2. In fact a simple linear least-squares fitting for the estimated values of $m_j(\cdot)$ led to the estimated slopes 0.2013, 0.2425, 0.2049 and 0.2552, for $j=1,2,3$ and 4, respectively, very close to the true values.

We also applied a *parametric bootstrap* method to test the null hypothesis of the auto-normal scheme (4.3): the bootstrap samples were generated from the auto-normal process with $\theta_1$ and $\theta_2$ in (4.3) estimated by the coding method (Besag [1]). Note that under the auto-normal scheme, $\mathrm{E}[\varepsilon(\mathbf{s})I\{\mathbf{X}(\mathbf{s})\in B\}]=0$ for any measurable $B\subset\mathbb{R}^4$, where $\varepsilon(\mathbf{s})$ is defined as the difference between $Y(\mathbf{s})$ and the right-hand side of (4.3), and $\mathbf{X}(\mathbf{s})$ consists of the four nearest neighbourhoods. This leads to the test statistic

$$T=\frac{1}{N}\sup_{1\leq k\leq N}\left|\sum_{j=1}^{N}\widehat{\varepsilon}(\mathbf{s}_j)I\{\mathbf{X}(\mathbf{s}_j)\leq\mathbf{X}(\mathbf{s}_k)\}\right|, \tag{4.4}$$

where $\mathbf{X}(\mathbf{s}_j)\leq\mathbf{X}(\mathbf{s}_k)$ is defined under the unilateral half-plane order (Whittle [25]), and

$$\widehat{\varepsilon}(\mathbf{s}_j)=Y(\mathbf{s}_j)-\widehat{\theta}_1\{Y(u_j-1,v_j)+Y(u_j+1,v_j)\}-\widehat{\theta}_2\{Y(u_j,v_j-1)+Y(u_j,v_j+1)\}.$$

Among the 500 replications, the proportions rejecting the linear null hypothesis are 10.8% at the level $\alpha=10\%$, and 4.4% at the level $\alpha=5\%$. To further assess the accuracy of the bootstrap approximation, we took the upper 10th and 5th percentiles for the empirical distribution of the 500 simulated values $T$ as the true critical values $t_\alpha$ for the test at the level $\alpha = 10\%$ and 5%, where Figure 3 displays the boxplots of the relative frequencies of the event $T^* > t_\alpha$ in the 200 bootstrap replications. This shows that most frequencies are clustered around $\alpha$ for both $\alpha = 10\%$ and 5%, indicating that the bootstrap approximation to the null distribution of $T$ is reasonably accurate.

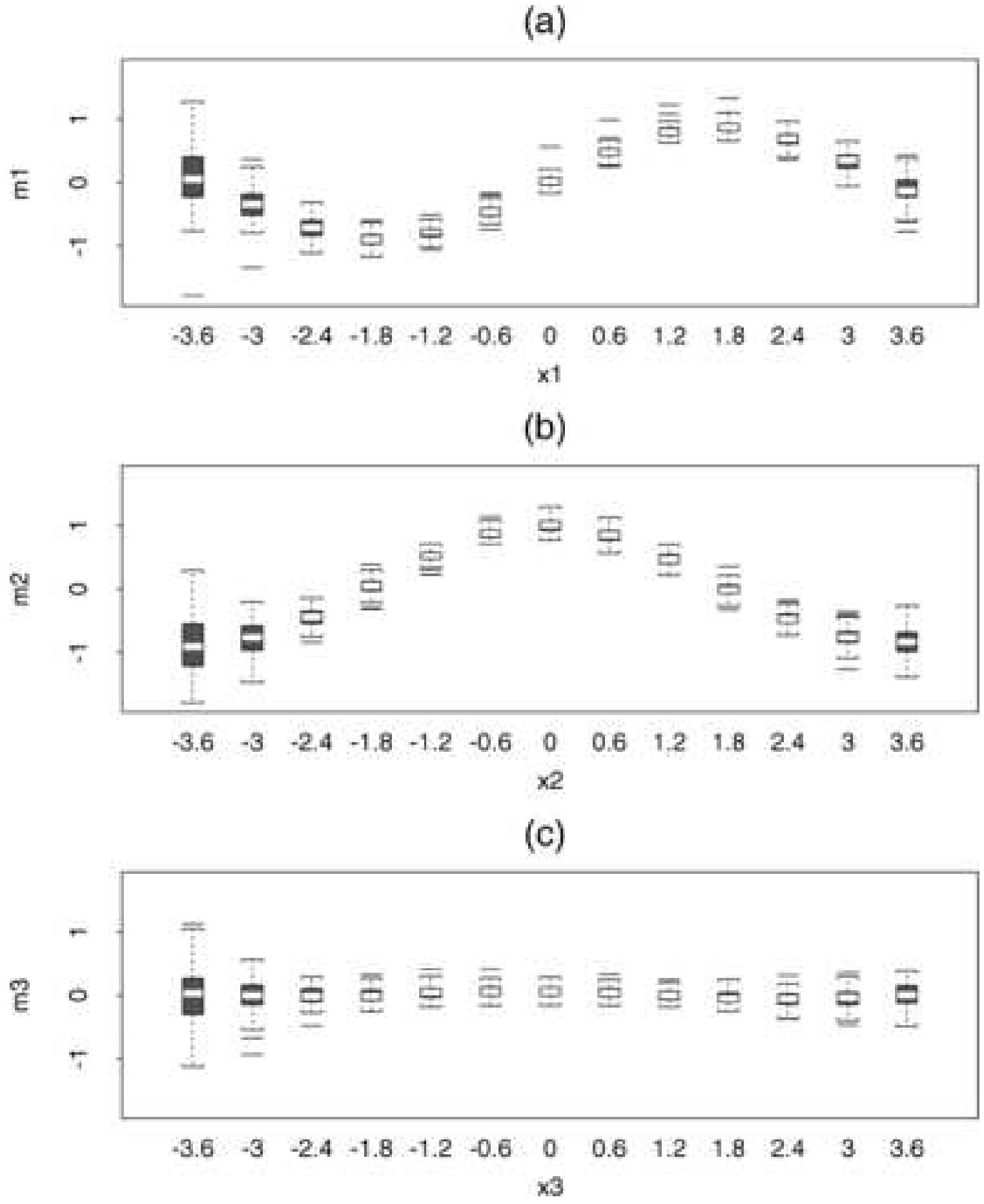

**Figure 1.** Boxplots of the estimators for (a) $m_1(x) = \sin(x) - \mathrm{E}[\sin\{Y(u,v)\}]$, (b) $m_2(x) = \cos(x) - \mathrm{E}[\cos\{Y(u,v)\}]$, and (c) $m_3(x) \equiv 0$ in Example 1.

## 4.2. A real data example

Figure 4 displays a magnetic resonance (MR) image of two test tubes filled with plastic straws of two different diameters embedded in gadolinium-doped agarose gel with tissue equivalent relaxation times. The straws test object was imaged with a T1-weighted SE pulse sequence on a Siemens Vision 1.5 T MR scanner using slice thickness of 4 mm and in-plane resolution of 0.6 mm × 0.6 mm. In the MR images of Figure 4, the more white a voxel is, the stronger the signal intensity. The black background region with very low intensity is air surrounding the two cylinders.

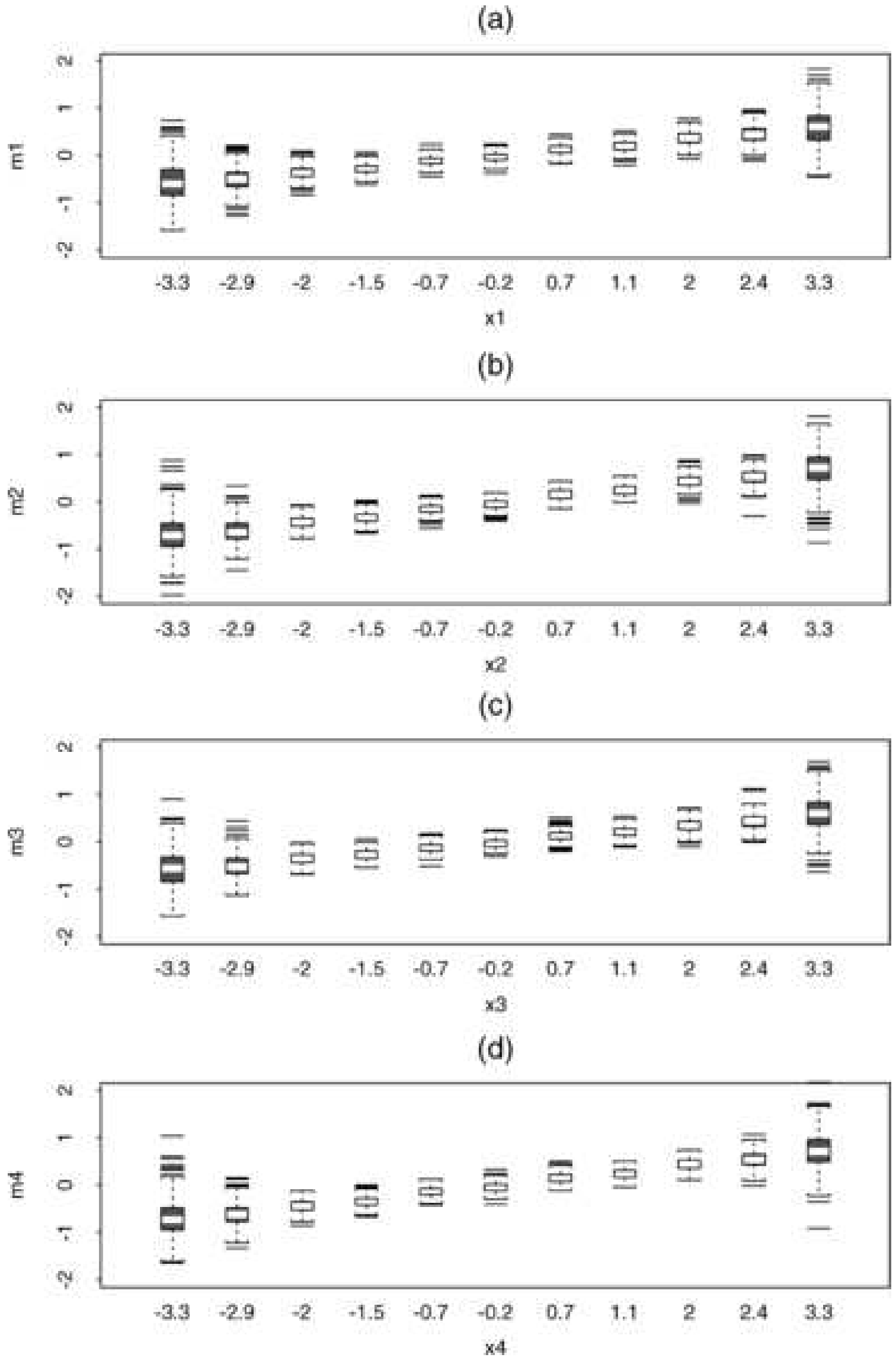

**Figure 2.** Boxplots of the estimators for (a) $m_1(x) = 0.2x$, (b) $m_2(x) = 0.25x$, (c) $m_3(x) = 0.2x$, and (d) $m_4(x) = 0.25x$ in Example 2.

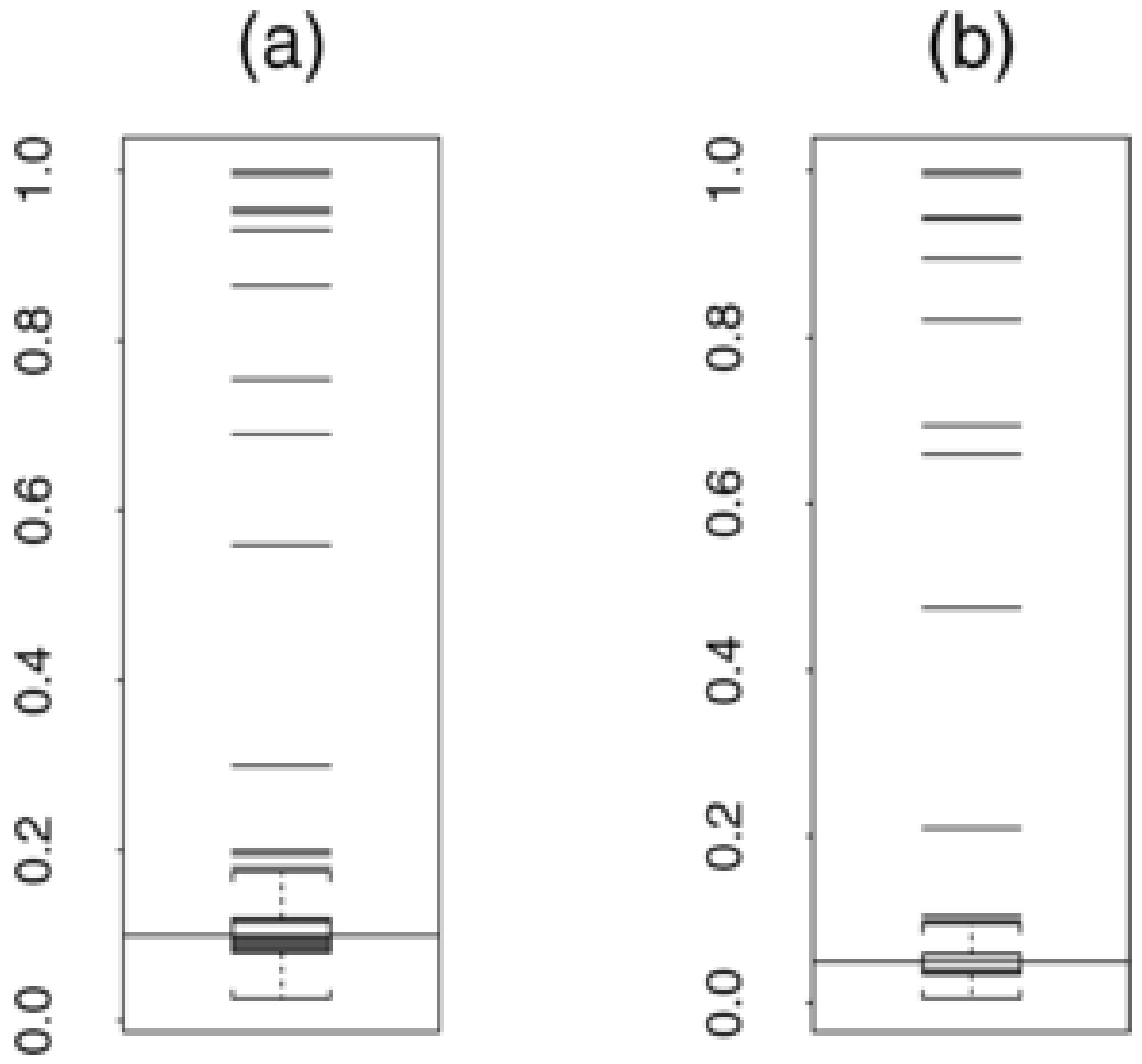

**Figure 3.** Boxplots of the relative frequencies of the event $\{T^* > t_\alpha\}$ for (a) $\alpha = 0.1$, and (b) $\alpha = 0.05$ in Example 2.

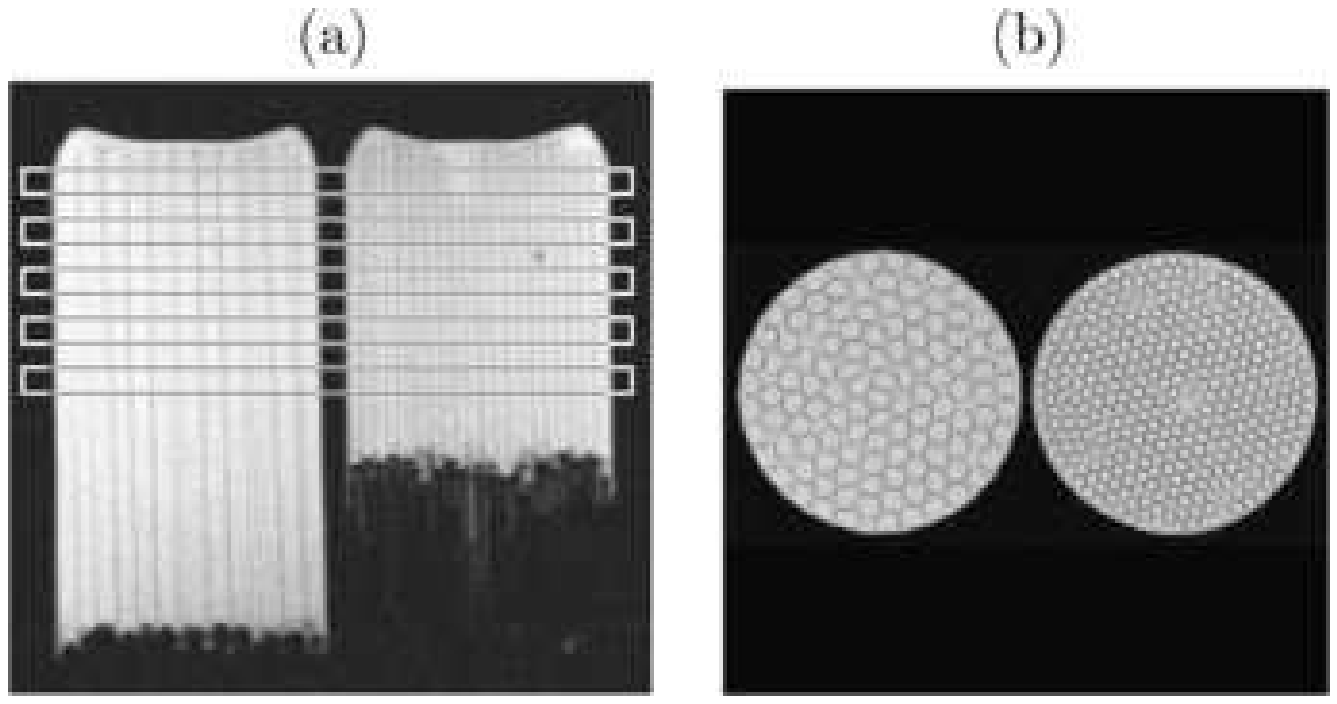

**Figure 4.** Modelling the straw data: MR images from the straws test object (a) depicting a longitudinal section with indication of the trans-axial slices, and (b) showing the upper trans-axial slice.

For our analysis, we chose two stationary-like subsets of image Figure 4(b), each of size 61×61. The subset images are plotted respectively in Figures 5(a)(i) and 6(a)(i). For each subset, we approximated the conditional expectation

$$E\{Y(u,v)|Y(u-1,v)=x_1, Y(u,v-1)=x_2, Y(u+1,v)=x_3, Y(u,v+1)=x_4\} \quad (4.5)$$

by an additive form

$$m_0 + m_1(x_1) + m_2(x_2) + m_3(x_3) + m_4(x_4),$$

where $m_1(x_1)$, $m_2(x_2)$, $m_3(x_3)$ and $m_4(x_4)$ represent the contributions from the nearest neighbourhood in the north, west, south and east direction, respectively. For comparison purposes, we also fitted the data using Besag's ([1]) first-order auto-normal scheme, assuming the conditional variances over different locations were the same. This leads to the assumption that the conditional expectation (4.5) is of the form

$$\alpha + \beta_1(x_1 - \alpha) + \beta_2(x_2 - \alpha) + \beta_1(x_3 - \alpha) + \beta_2(x_4 - \alpha). \tag{4.6}$$

The coefficients $\beta_j$ and $\alpha$ were estimated using Besag's coding method. The estimated additive functions $\widetilde{m}_j(x)$, together with the fitted straight lines $\widehat{\beta}_j(x - \widehat{\alpha})$, are plotted in Figure 5 for the large-diameter straws and in Figure 6 for the small-diameter straws. As in Example 2 above, we also applied the parametric bootstrap based on statistic (4.4) for testing the auto-normal null hypothesis for those two subsets (with 200 bootstrap replications), leading to $p$-values less than 0.005. This indicates that the first order auto-normal scheme is inadequate for both the data sets. The histograms presented in Figures 5(a)(ii) and 6(a)(ii) indicate bimodal marginal density functions; the lower-intensity mode corresponds to voxels at the boundary between straws, and the higher-intensity mode to voxels in the interior of the straws. The pointwise confidence intervals for $m_j(\cdot)$ were obtained using the standard normal $\varepsilon(\mathbf{s})$ in (3.13) with 100 wild bootstrap replications; see (3.14).

The plots of the $m_j(\cdot)$ functions as a rough approximation suggest isotropy, which is natural given the set-up of the straws. Both the plots and bootstrap test point to nonlinearity. Note that the bends at the ends of the curves cannot be attributed to boundary effects, as substantial number of voxels fall in the end regions; see the histograms in Figures 5 and 6.

A possible explanation for the bends is as follows: for both the large- and small-diameter straws there is a positive correlation among intensity values in the middle (the $m_j(\cdot)$ curve has a positive slope). These intensity values correspond to voxels in the center of the straws, for which it is seen (from Figures 5(a)(i) and 6(a)(i)) that there is a positive spatial autocorrelation. For the small-diameter straws there is local negative correlation at both ends of the curves. Looking at the image (Figure 6(a)(i)) it is seen that the lowest-intensity values (darkest) voxels are at the boundaries as expected. But it is also seen that there are voxels of very high intensity (very white) close to the boundary. Similarly, the voxels of highest intensity are often found close to the boundary and have low-intensity voxels in their neighbourhood, resulting in the local negative correlation for extreme intensity values. For large-diameter straws with low-intensity voxels at the boundary, the same pattern occurs but not to the same extent; see Figure 5(a)(i). The picture for voxels of high intensity is less clear, as some of these are located close to the boundary surrounded by low-intensity voxels, others close to the centre with high-intensity neighbours. There is no clear dependence pattern for high values, which is echoed by the flatness of the $m_j(\cdot)$ plots for high intensities.

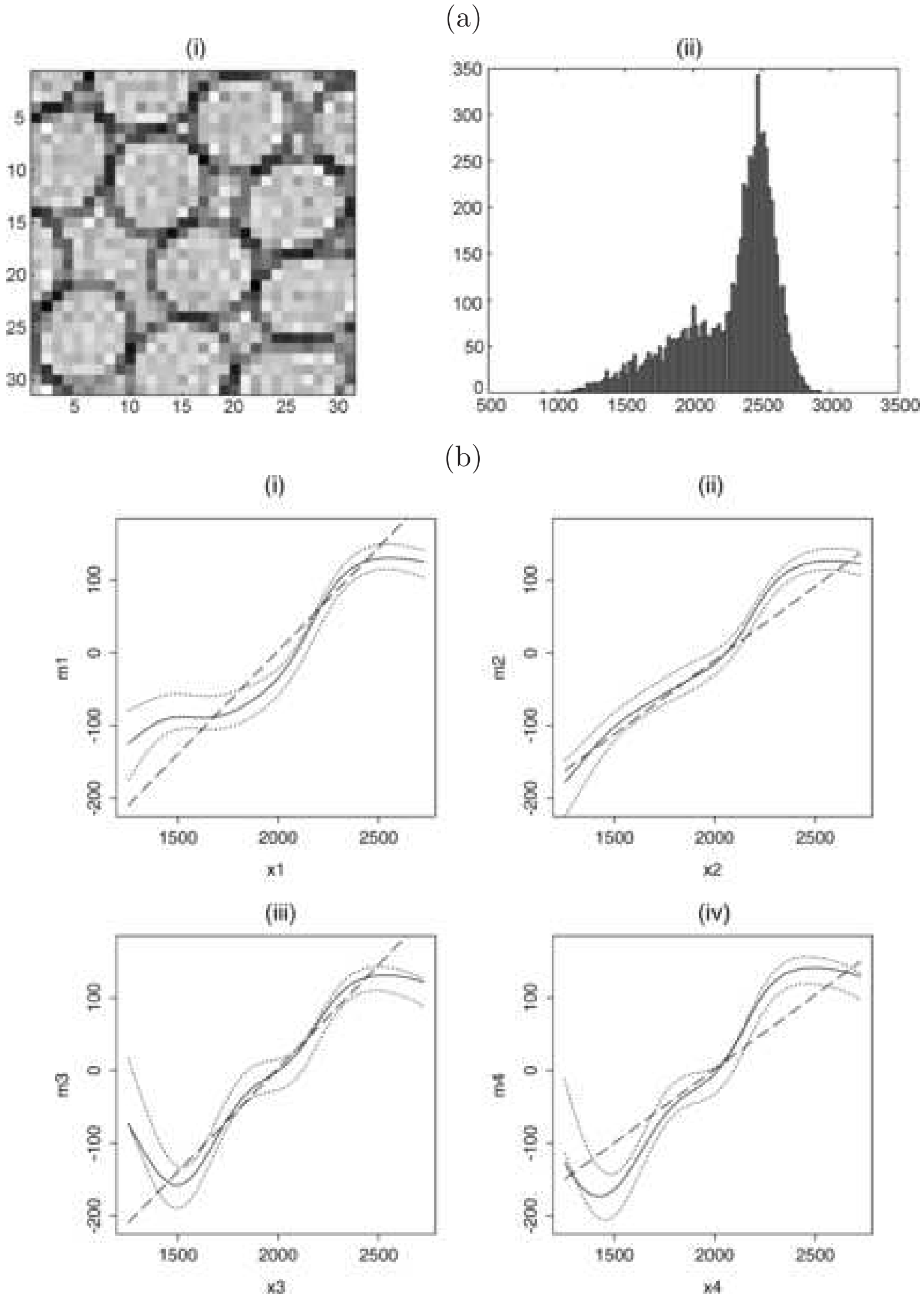

**Figure 5.** Modelling subset I of straw data. (a) (i) A subregion of interest ($61 \times 61$ window, of which a $31 \times 31$ detail is shown) from the bundle of large diameter straws in Figure 4(b); (ii) the corresponding signal intensity histogram. Pixel signal intensity in a.u. (12-bit range). (b) Additive estimators (solid lines), the boundaries of pointwise confidence intervals (dotted lines), and auto-normal scheme estimators (dashed lines) for (i) $m_1(x)$, (ii) $m_2(x)$, (iii) $m_3(x)$ and (iv) $m_4(x)$ and $m_0 = 2278.615$.

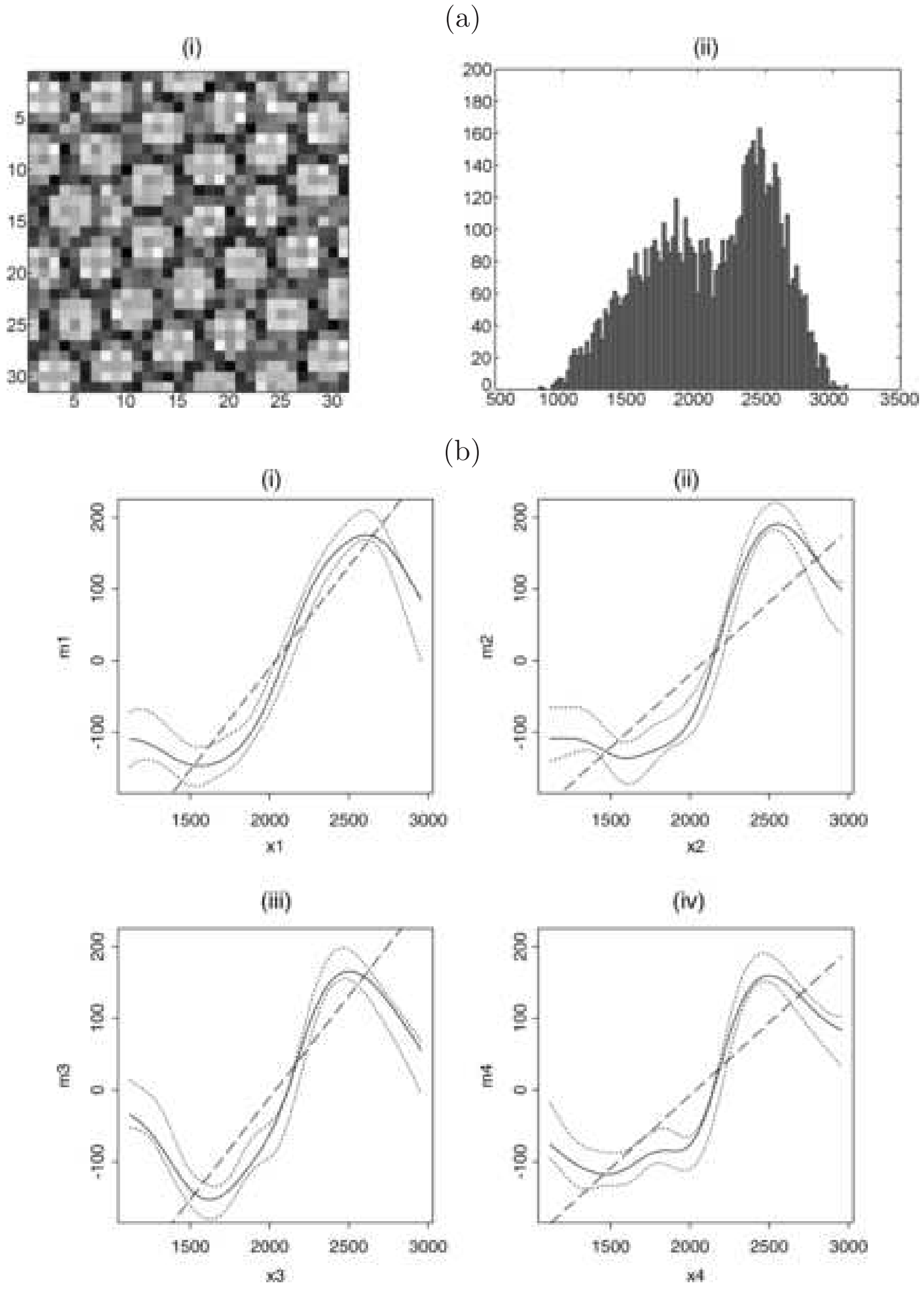

**Figure 6.** Modelling subset II of straw data. (a) (i) A subregion of interest ($61 \times 61$ window, of which a $31 \times 31$ detail is shown) from the bundle of small diameter straws in Figure 4(b); (ii) the corresponding signal intensity histogram. Pixel signal intensity in a.u. (12-bit range). (b) Additive estimators (solid lines), boundaries of pointwise confidence intervals (dotted lines) and auto-normal scheme estimators (dashed lines) for (i) $m_1(x)$, (ii) $m_2(x)$, (iii) $m_3(x)$ and (iv) $m_4(x)$ and $m_0 = 2089.465$.

Overall the nonparametric additive approximations for conditional means describe the local correlation structure of the straws quite well, whereas the auto-normal models fail to do so since they only reproduce the dominating positive spatial autocorrelation in the interior of the straws.

## 5. Discussion

Observations taken on irregular grids often occur in practical spatial problems. We outline below some tentative ideas to extend the method proposed in this paper to handle irregularly spaced data. Our basic assumption is that the observations $\{Y(\mathbf{s}_j), j=1,\dots,N\}$ (after detrending appropriately) are taken over an irregular grid from a strictly stationary process $Y(\mathbf{s})$ with index $\mathbf{s}$ varying continuously on $\mathbb{R}^2$. Our goal is to estimate the best additive approximation, in the sense of (2.4) and (2.5), for the conditional expectation at a fixed location given its $d$ neighbourhood observations, for small $d$ such as 3 or 4. Without loss of generality, we may assume that the location is at $\mathbf{0}=(0,0)^\top$, and the $d$ neighbourhood locations are $\mathbf{i}_1,\dots,\mathbf{i}_d$. Put $\mathbf{X}(\mathbf{0})=\{Y(\mathbf{i}_1),\dots,Y(\mathbf{i}_d)\}^\top$. Our task is now to estimate the best additive approximation for

$$m(\mathbf{x})=\mathrm{E}\{Y(\mathbf{0})|\mathbf{X}(\mathbf{0})=\mathbf{x}\}.$$

First, for each location $\mathbf{s}_k$ at which we have an observation $Y(\mathbf{s}_k)$, we define its $d$ neighbourhoods selected among the other $N-1$ observations by minimizing

$$\sum_{j=1}^{d}\|\mathbf{i}_j-(\mathbf{s}^j-\mathbf{s}_k)\|,$$

where the minimization is taken over $\mathbf{s}^j\in\{\mathbf{s}_1,\dots,\mathbf{s}_N\}$, $\mathbf{s}^j\neq\mathbf{s}_k$, and $\mathbf{s}^1,\dots,\mathbf{s}^d$ are all different from each other. Let $(\mathbf{s}_{k1},\dots,\mathbf{s}_{kd})$ be the minimizer. Then $\mathbf{X}(\mathbf{s}_k)=\{Y(\mathbf{s}_{k1}),\dots,Y(\mathbf{s}_{kd})\}^\top$ are the $d$ neighbourhood observations of $Y(\mathbf{s}_k)$ as far as our task is concerned. Put

$$\lambda_k=\sum_{j=1}^{d}\|\mathbf{i}_j-(\mathbf{s}_{jk}-\mathbf{s}_k)\|,$$

which measures the discrepancy between the pattern of $(\mathbf{s}_{1k},\dots,\mathbf{s}_{dk})$ in relation to $\mathbf{s}_k$ and that of $(\mathbf{i}_1,\dots,\mathbf{i}_d)$ in relation to $\mathbf{0}$. It is easy to see that $\lambda_k=0$ if and only if $(\mathbf{s}_k,\mathbf{s}_{1k},\dots,\mathbf{s}_{dk})$ is merely a shift (without rotating) of $(\mathbf{0},\mathbf{i}_1,\dots,\mathbf{i}_d)$ in $\mathbb{R}^2$. The larger $\lambda_k$ is, the larger the discrepancy is between the two patterns. As far as the estimation for $m(\mathbf{x})$ is concerned, we should not treat all $\{Y(\mathbf{s}_k),\mathbf{X}(\mathbf{s}_k)\}$ equally, as we did for regularly spaced data. Instead we give more weight to the observations $\{Y(\mathbf{s}_k),\mathbf{X}(\mathbf{s}_k)\}$ with smaller values of $\lambda_k$. By taking this into account, an argument similar to that in Section 2.2 leads to the backfitting estimation (2.12) in which, however, now

$$\widehat{f}_j(x_j)=\sum_{k=1}^{N}w_kK_h\{x_j-Y(\mathbf{s}_{jk})\},\qquad \widehat{m}_j(x_j)=\frac{1}{\widehat{f}_j(x_j)}\sum_{k=1}^{N}w_kY(\mathbf{s}_k)K_h\{x_j-Y(\mathbf{s}_{jk})\},$$

$$\widehat{f}_{j\ell}(x_j, x_\ell) = \sum_{k=1}^{N} w_k K_h\{x_j - Y(\mathbf{s}_{jk})\} K_h\{x_\ell - Y(\mathbf{s}_{\ell k})\},$$

where the weight function $w_k = W(\lambda_k/b)/\sum_{j=1}^{N} W(\lambda_j/b)$, $W(\cdot)$ is a kernel function, and $b > 0$ is a bandwidth. The associated issues on inference, computation and asymptotic properties are subject to further investigation. An alternative that could be explored is to replace $\mathbf{i}_1, \dots, \mathbf{i}_d$ by an average set of $d$ neighbourhood points $\mathbf{i}'_1, \dots, \mathbf{i}'_d$, where $\mathbf{i}'_j = N^{-1}\sum_{k=1}^{N}(\mathbf{s}'_{kj} - \mathbf{s}_k)$, $j = 1, \dots, d$, in which $\mathbf{s}'_{kj}$ is the $j$th nearest neighbour of $\mathbf{s}_k$.

Finally, we note that for observations taken irregularly over space and regularly over time, the method proposed in Section 2.2 may be applied directly if we only use the data taken at the fixed location but over different times in the estimation. Technically this reduces to a problem of multivariate time series modelling. However, there is an added advantage: the inference does not rely on the assumption of the stationarity over space.

# Appendix

## A.1. Proof of Theorem 1

We only need to justify conditions (A1)–(A3) in Mammen *et al.* [15]. The required result follows from their Theorem 1 immediately.

Condition (A1) requires that for $j \neq k$,

$$\int \frac{p_{jk}^2(x_j, x_k)}{p_k(x_k) p_j(x_j)} \, \mathrm{d}x_j \, \mathrm{d}x_k = \int_{A_{ij}} \frac{p_{jk}^2(x_j, x_k)}{p_k(x_k) p_j(x_j)} \, \mathrm{d}x_j \, \mathrm{d}x_k < \infty,$$

which is guaranteed by (C5). By Theorem 2 of Yao [26],

$$\sup_{x_j \in A_j} |\widehat{p}_j(x_j) - p_j(x_j)| = O_p\left\{h^2 + \left(\frac{\log N}{Nh}\right)^{1/2}\right\}. \tag{A.1}$$

Similarly, we may show that

$$\sup_{(x_j, x_k) \in A_{jk}} |\widehat{p}_{jk}(x_j, x_k) - p_{jk}(x_j, x_k)| = O_p\left\{h^2 + \left(\frac{\log N}{Nh^2}\right)^{1/2}\right\}. \tag{A.2}$$

Furthermore, it is easy to see from Theorem 3 below that

$$\sup_{x_j \in A_j} |\widehat{m}_j(x_j) - \mathrm{E}\widehat{m}_j(x_j)| = O_p(1), \qquad \sup_{x_j \in A_j} |\mathrm{E}\widehat{m}_j(x_j)| \leq C. \tag{A.3}$$

Note that $f_{jk}(x_j, x_k) = f_{kj}(x_k, x_j)$. Condition (C5) implies that $(x_j, x_k) \in A_{jk}$ if and only if $(x_k, x_j) \in A_{kj}$ for any $j \neq k$. This, together with (A.1)–(A.3), entails conditions (A2) and (A3) of Mammen *et al.* [15].

### A.2. Proof of Theorem 2

Let $e(\mathbf{s}) = Y(\mathbf{s}) - m\{\mathbf{X}(\mathbf{s})\}$, and

$$\widehat{m}_j^a(x_j) = \frac{I(x_j \in A_j)}{\widehat{p}(x_j)\sum_{\ell=1}^N I\{Y(\mathbf{s}_\ell - \mathbf{i}_j) \in A_j\}} \times \sum_{\ell=1}^N [e(\mathbf{s}_\ell) + m\{\mathbf{X}(\mathbf{s}_\ell)\} - m^o\{\mathbf{X}(\mathbf{s}_\ell)\}] K_h\{x_j - Y(\mathbf{s}_\ell - \mathbf{i}_j)\}, \tag{A.4}$$

$$\widehat{m}_j^b(x_j) = \frac{I(x_j \in A_j)}{\widehat{p}(x_j)\sum_{\ell=1}^N I\{Y(\mathbf{s}_\ell - \mathbf{i}_j) \in A_j\}} \sum_{\ell=1}^N m^o\{\mathbf{X}(\mathbf{s}_\ell)\} K_h\{x_j - Y(\mathbf{s}_\ell - \mathbf{i}_j)\}.$$

Then $\widehat{m}_j(x_j) = \widehat{m}_j^a(x_j) + \widehat{m}_j^b(x_j)$. Let $\widetilde{m}_j^a(x_j)$ and $\widetilde{m}_j^b(x_j)$ be defined by (3.6) with $\widehat{m}_j(x_j)$ replaced by $\widehat{m}_j^a(x_j)$ and $\widehat{m}_j^b(x_j)$, respectively. We first introduce two lemmas.

**Lemma 1.** *Under conditions* (C1)–(C5), *for any* $j \neq k$,

$$\sum_{x_k \in A_k} \left| \int \frac{\widehat{p}_{jk}(x_j, x_k)}{\widehat{p}_k(x_k)} \widehat{m}_j^a(x_j)\, \mathrm{d}x_j \right| = o_p(h^2)$$

*and*

$$\int p_k(x_k)\, \mathrm{d}x_k \left\{ \int \frac{\widehat{p}_{jk}(x_j, x_k)}{\widehat{p}_k(x_k)} \widehat{m}_j^a(x_j)\, \mathrm{d}x_j \right\}^2 = o_p(h^4).$$

**Lemma 2.** *Under conditions* (C1)–(C5), $\int m_j^o(x_j)\widehat{p}_j(x_j)\, \mathrm{d}x_j = o_p(h^2)$, *and*

$$\sup_{x_k \in A_k} |\widehat{m}_j^b(x_j) - \widehat{\mu}_j(x_j)| = o_p(h^2), \qquad \int |\widehat{m}_j^b(x_j) - \widehat{\mu}_j(x_j)|^2 p_j(x_j)\, \mathrm{d}x_j = o_p(h^4).$$

Based on (3.11), Lemma 1 may be proved in the same manner as (A6) in Appendix A of Mammen *et al.* [15]. The proof of Lemma 2 is similar to the proofs of (114), (112) and (113) in Mammen *et al.* [15].

We now sketch the proof of Theorem 2. Note that $\mathbf{x} = (x_1, \ldots, x_d)^\top$ is an inner point of $A$. Lemma 1 implies condition (A6) of Mammen *et al.* [15] with $\Delta_N = h^2$. By Theorem 3 below, condition (A9) of Mammen *et al.* [15] also holds. By Theorem 3 of Mammen *et al.* [15], condition (A7) in their paper also holds. It may be proved that

$$\int \widehat{m}_j^a(x)\widehat{p}_j(x)\, \mathrm{d}x \bigg/ \int \widehat{p}_j(y)\, \mathrm{d}y = O_p(N^{-1/2}) = o_p(h^2).$$

Now it follows from Theorems 2 and 3 of Mammen *et al.* [15] that

$$\widetilde{m}_j(x_j) = \widehat{m}_j^a(x_j) + h^2\beta_j(x_j) + o_p(h^2). \tag{A.5}$$

Note that

$$\mathrm{E}[e(\mathbf{s}_\ell)K_h\{x_j - Y(\mathbf{s}_\ell - \mathbf{i}_j)\}] = \mathrm{E}[\mathrm{E}\{e(\mathbf{s}_\ell)|\mathbf{X}(\mathbf{s}_\ell)\}K_h\{x_j - Y(\mathbf{s}_\ell - \mathbf{i}_j)\}] = 0$$

and

$$\begin{aligned} &E[\{m\{\mathbf{X}(\mathbf{s}_\ell)\} - m^o\{\mathbf{X}(\mathbf{s}_\ell)\}\}K_h\{z - Y(\mathbf{s}_\ell - \mathbf{i}_j)\}] \\ &\quad = \int K_h(z - x_j)\left[\int \{m(\mathbf{x}) - m^o(\mathbf{x})\}p(\mathbf{x})\prod_{k\neq j} \mathrm{d}x_k\right] \mathrm{d}x_j = 0. \end{aligned} \tag{A.6}$$

The last equality in the above expression follows from (3.11). Now the required central limit theorem follows from (A.5), (A.4), (A.6) and the theorem in Bolthausen [2].

## A.3. Uniform convergence rates for regression estimation

First, we introduce some notation. Let $\{(Y(\mathbf{s}_j), X(\mathbf{s}_j)), 1 \le j \le N\}$ be observations from a two-dimensional strictly stationary spatial process with $\{\mathbf{s}_1, \ldots, \mathbf{s}_N\}$ given as in (2.7). Let $f(\cdot)$ be the density function of $X(\mathbf{s})$ and $m(x) = \mathrm{E}\{Y(\mathbf{s})|X(\mathbf{s}) = x\}$. We define the Nadaraya–Watson estimator $\widehat{m}(x) = \widehat{r}(x)/\widehat{f}(x)$ with

$$\widehat{f}(x) = \frac{1}{N}\sum_{j=1}^{N} K_h\{x - X(\mathbf{s}_j)\}, \qquad \widehat{r}(x) = \frac{1}{N}\sum_{j=1}^{N} Y(\mathbf{s}_j)K_h\{x - X(\mathbf{s}_j)\}.$$

We introduce some regularity conditions.

- (C1′) $m(\cdot)$ has continuous first derivative, $f(\cdot)$ has continuous second derivative, and the joint density function of $X(\mathbf{s})$ and $X(\mathbf{s}+\mathbf{i})$ is bounded by a constant independent of $\mathbf{i}$. Furthermore, (3.2) holds.
- (C4′) $\alpha(k; k', j) \le Ck^{-\beta}$ for any $k, j$ and $k' = O(k^2)$, where $\alpha$ is defined as in (3.1) with $\mathbf{X}(\mathbf{s})$ replaced by $X(\mathbf{s})$.

**Theorem 3.** *Let $A$ be a compact set contained in the support of $f(\cdot)$. Under conditions* (C1′), (C2), (C3) *and* (C4′),

$$\sup_{x\in A} |\widehat{m}(x) - m(x)| = O_p\left\{\left(\frac{\log N}{Nh}\right)^{1/2} + h^2\right\} \tag{A.7}$$

*and*

$$\sup_{x\in A} |\widehat{m}(x) - \mathrm{E}\widehat{m}(x)| = O_p\left\{\left(\frac{\log N}{Nh^2}\right)^{1/2} + h^2\right\}. \tag{A.8}$$

**Proof.** We first prove (A.7). Let $r(x)=m(x)f(x)$. For any $a_N>0$ and $\varepsilon>0$,

$$
\begin{aligned}
&P\Big\{\sup_{x\in A} a_N|\widehat{m}(x)-m(x)|>\varepsilon\Big\}\\
&\quad\le P\bigg\{\frac{a_N}{\inf_{y\in A}\widehat{f}(y)}\bigg(\sup_{x\in A}|\widehat{r}(x)-r(x)|+\max_{z\in A} m(z)\sup_{x\in A}|\widehat{f}(x)-f(x)|\bigg)>\varepsilon\bigg\}\\
&\quad\le P\Big\{C_1a_N\sup_{x\in A}|\widehat{r}(x)-r(x)|+C_2a_N\sup_{x\in A}|\widehat{f}(x)-f(x)|>\varepsilon\Big\}\\
&\qquad+P\Big\{\inf_{x\in A}|\widehat{f}(x)-f(x)|>\tau\Big\},
\end{aligned}
$$

where $C_1, C_2, \tau>0$ are constants. It follows from Theorem 2 of Yao [26] that the second term of the right-hand side of the above expression may be arbitrarily small for all sufficiently large $N$. Then (A.7) follows from Theorem 2 of Yao [26] and

$$
\sup_{x\in A}|\widehat{r}(x)-r(x)|=O_p\bigg\{\bigg(\frac{\log N}{Nh}\bigg)^{1/2}+h^2\bigg\}, \tag{A.9}
$$

which will now be established.

Partition $A$ into $L$ subintervals $\{I_j\}$ of equal length. Let $x_j$ be the centre of $I_j$. Since

$$
|\widehat{r}(x)-\widehat{r}(x')|\le\frac{1}{N}\sum_{j=1}^{N}|Y(\mathbf{s}_j)||K_h\{X(\mathbf{s}_j)-x\}-K_h\{X(\mathbf{s}_j)-x'\}|\le\frac{1}{N}\sum_{j=1}^{N}|Y(\mathbf{s}_j)|\frac{C}{h}|x-x'|,
$$

we have that $|\mathrm{E}\widehat{r}(x)-\mathrm{E}\widehat{r}(x')|\le Ch^{-1}|x-x'|$. Hence

$$
\sup_{x\in A}|\widehat{r}(x)-\mathrm{E}\widehat{r}(x)|\le\max_{1\le j\le L}|\widehat{r}(x_j)-\mathrm{E}\widehat{r}(x_j)|+O_p\bigg(\frac{1}{Lh}\bigg). \tag{A.10}
$$

For large $M>0$, define

$$
\widehat{r}_1(x)=\frac{1}{N}\sum_{j=1}^{N}Y(\mathbf{s}_j)I\{|Y(\mathbf{s}_j)|\le M\}K_h\{X(\mathbf{s}_j)-x\},
$$

$$
\widehat{r}_2(x)=\frac{1}{N}\sum_{j=1}^{N}Y(\mathbf{s}_j)I\{|Y(\mathbf{s}_j)|>M\}K_h\{X(\mathbf{s}_j)-x\}
$$

and

$$
r_1(x)=\mathrm{E}[Y(\mathbf{s}_j)I\{|Y(\mathbf{s}_j)|\le M\}|X(\mathbf{s}_j)=x],\qquad r_2(x)=\mathrm{E}[Y(\mathbf{s}_j)I\{|Y(\mathbf{s}_j)|>M\}|X(\mathbf{s}_j)=x].
$$

Then $\widehat{r}(x)=\widehat{r}_1(x)+\widehat{r}_2(x)$ and $r(x)=r_1(x)+r_2(x)$. Since $|K_h(\cdot)|\le Ch^{-1}$, it follows from the second inequality in Theorem 1 of Yao [26] that

$$P\{|\widehat{r}_1(x)-\mathrm{E}\widehat{r}_1(x)|>\varepsilon\}\le 8\exp\left\{-\frac{\varepsilon^2q^2}{8\nu(q)^2}\right\}+44\left(1+\frac{4CM}{\varepsilon h}\right)^{1/2}q^2\alpha([p_1]\wedge[p_2];[p_1p_2],N),$$

where $q=[(\varepsilon M)^{1/2}(N_1\wedge N_2)]$, $p_i=N_i/(2q)$ and

$$\nu(q)^2\le\frac{32q^4}{N^2}\frac{Cp_1p_2}{h}+\frac{CM\varepsilon}{h}=\frac{C_1}{p_1p_2h}+\frac{CM\varepsilon}{h}\le\frac{C_2M\varepsilon}{h},$$

where $C,C_1,C_2>0$ are constant. The first inequality in the above expression can be verified similarly to the variance expression in Proposition 1 of Yao [26], and the second inequality is obvious by setting $M=\log N$ and $\varepsilon^2=8aC\log N/(N_1\wedge N_2)^2h$ for some constant $a>0$. Now

$$\exp\left\{-\frac{\varepsilon^2q^2}{8\nu(q)^2}\right\}\le\exp\left\{-\frac{\varepsilon^2(N_1\wedge N_2)^2h}{8C}\right\}=N^{-a}.$$

On the other hand, condition (C4′) entails that

$$\begin{aligned}\left(\frac{M}{\varepsilon h}\right)^{1/2}q^2\alpha(p_1\wedge p_2;[p_1p_2],N)&\le C\left(\frac{M}{\varepsilon h}\right)^{1/2}\varepsilon MN(\varepsilon M)^{\beta/2}\\&=CM^{(\beta+3)/2}Nh^{-1/2}\varepsilon^{(\beta+1)/2}\\&=O\{N^{-\beta/4+3/4}h^{-\beta/4-3/4}(\log N)^{3\beta/4+7/4}\}.\end{aligned}$$

Let $L=[(N/h)^{1/2}]$. Hence,

$$\begin{aligned}&P\left\{\max_{1\le j\le L}|\widehat{r}_1(x_j)-\mathrm{E}\widehat{r}_1(x_j)|>\varepsilon\right\}\\&\quad\le L\{N^{-a}+N^{-\beta/4+3/4}h^{-\beta/4-3/4}(\log N)^{3\beta/4+7/4}\}\to 0,\end{aligned}\tag{A.11}$$

see condition (3.3). On the other hand,

$$P\{|\widehat{r}_2(x)-\mathrm{E}\widehat{r}_2(x)|>\varepsilon\}\le NP\{|Y(\mathbf{s})|>M\}\le N\mathrm{e}^{-\lambda M}\mathrm{E}\{\mathrm{e}^{\lambda|Y(\mathbf{s})|}\}=O(N^{-\lambda+1}),$$

where $\lambda>0$ is a constant. Hence

$$P\left\{\max_{1\le j\le L}|\widehat{r}_2(x)-\mathrm{E}\widehat{r}_2(x)|>\varepsilon\right\}\le O(LN^{-\lambda+1})\to 0;$$

see (3.2). Combining this with (A.11) and (A.10), we have

$$P\left\{\sup_{x\in A}|\widehat{r}(x)-\mathrm{E}\widehat{r}(x)|>\varepsilon\right\}=O_p\left\{\left(\frac{\log N}{Nh}\right)^{1/2}\right\}.$$

Now (A.9) follows from this and the fact that $\sup_{x\in A}|\mathrm{E}\widehat{r}(x) - r(x)| = O(h^2)$, which may be established via simple algebraic manipulation.

To prove (A.8), it follows from (A.7) and Theorem 2 of Yao [26] that for any $\varepsilon > 0$ there exists a $\tau > 0$ such that

$$\sup_{x\in A}|\mathrm{E}\widehat{m}(x) - \mathrm{E}\{\widehat{m}(x)I(|\widehat{f}(x) - f(x)| < h^2\tau)\}| < \varepsilon.$$

Now, for $x \in A$,

$$\begin{aligned}
&\sup_{x\in A}|\mathrm{E}\{\widehat{m}(x)I(|\widehat{f}(x) - f(x)| < h^2\tau)\} - m(x)| \\
&\quad \le \sup_{x\in A}|\mathrm{E}\{\widehat{r}(x)/f(x)I(|\widehat{f}(x) - f(x)| < h^2\tau)\} - m(x)| + Ch^2 \\
&\quad \le \sup_{x\in A}|\mathrm{E}\{\widehat{r}(x)/f(x)\} - m(x)| + \varepsilon + C_1h^2 \le \varepsilon + C_2h^2.
\end{aligned}$$

holds uniformly. Hence (A.8) follows from (A.7). □

## Acknowledgements

This project was partially supported by a Leverhulme Trust research grant. Lu's work was also partially supported by the National Natural Science Foundation of China (70271003, 70221001). The authors thank Dr. Michael Bock for providing the MR straws data as part of the EC COST B11 project 'Quantification of Magnetic Resonance Image Texture', Professor Enno Mammen for helpful discussion on the asymptotic theory in Section 4, and two referees for very helpful comments and suggestions.